\documentclass[11pt]{article}

\usepackage{amsfonts}
\usepackage [dvips]{graphics}
\usepackage[centertags]{amsmath}
\usepackage{amsfonts}
\usepackage{amssymb}
\usepackage{amsthm}
\usepackage{newlfont}
\usepackage{mathrsfs}
\usepackage [dvips]{graphics}
\usepackage[centertags]{amsmath}
\usepackage{amsfonts}
\usepackage{amssymb}
\usepackage{amsthm}
\usepackage{newlfont}
\usepackage{hyperref}

\newlength{\defbaselineskip}
\newcommand{\setlinespacing}[1]%
           {\setlength{\baselineskip}{#1 \defbaselineskip}}

\theoremstyle{plain}
\newtheorem{thm}{Theorem}[section]

\newtheorem{lem}[thm]{Lemma}

\theoremstyle{definition}

\makeatletter\@addtoreset{equation}{section} \makeatother

\begin{document}

\title{ Optimal control problem
of fully  coupled forward-backward stochastic
 systems with Poisson jumps under partial information    \footnotemark[1] }
\author{Qingxin MENG\footnotemark[2]}  
\date{}

\footnotetext[1]{This work was partially supported by Basic Research
Program of China(Grant No. 2007CB814904) and the National Natural
Science Foundation of China (Grant No.10325101 )and the Natural
Science Foundation of Zhejiang Province
(Grant No. Y605478, Y606667).}
\footnotetext[2]{Institute of Mathematics, Fudan University,
Shanghai 200433, China, Email: 071018034@fudan.edu.cn}

\maketitle

\begin{abstract}  In this paper, we study a class of  stochastic
optimal control problem with jumps under partial information. More
precisely, the controlled systems are  described by a fully coupled
nonlinear multi- dimensional forward-backward stochastic
differential equation driven by a Poisson random measure  and an
independent multi-dimensional Brownian motion, and all admissible
control processes are required to be adapted to a given
subfiltration of the filtration generated by the underlying Poisson
random measure and Brownian motion. For this type of partial
information stochastic optimal control problem, we give a necessary
and sufficient maximum principle. All the coefficients appearing in
the systems are allowed to depend on the control variables and the
control domain is convex.

\end{abstract}


\maketitle

\section{ Introduction } \vskip 0.5cm

  In recent years, there have been growing interests on stochastic
optimal control problems under partial information, partly due to
the applications in mathematical finance. For the partial
information  optimal control problem, the objective is to find an
optimal control for which the controller has less information than
the complete information filtration. In particular, sometimes an
economic model in which there are information gaps among economic
agents can be formulated as a partial information optimal control
problem (see {\O}ksendal\cite{OKSN} , Kohlmann and Xiong
\cite{KoXi}).

There are two important approaches to the general stochastic optimal
control problem. One is the Bellman dynamic programming principle,
which results in the Hamilton-Jacobi-Bellman equation. Another
important approach is the stochastic maximum principle by which  one
necessary or one sufficient condition of optimality can be obtained
by duality theory. For detailed accounts of the approaches for
complete information stochastic optimal control problem of the
forward system, see the books \cite{ZhYo} and the references
therein.

  Recently, Baghery and
${\O}$ksendal \cite{BaOk} established a maximum principle of forward
systems with jumps for under partial information. In \cite{BaOk},
the authors point out that because of the general nature of the
partial information filtration, dynamic programming and
Hamilton-Jacobi-Bellman equation can not be used to solve the
corresponding stochastic optimal control problem.

Backward stochastic differential equations coupled with forward
stochastic differential equations are called forward- backward
stochastic differential equations (FBSDEs). Forward-Backward
stochastic systems are not only encountered in  stochastic optimal
control problem when applying the stochastic maximum principle but
also used in mathematical finance (see Antonelli \cite{Anto}, Duffie
and Epstein \cite{DuEp}, El Karoui, Peng and Quenez \cite{ElPe} for
example). It now becomes more clear that certain important problems
in mathematical economics and mathematical finance, especially in
the optimization problem, can be formulated to be FBSDEs. In 2009,
Meng \cite{Meng}studied the partial stochastic optimal control
problem of continuous fully coupled forward-backward stochastic
systems driven by a Brownian motion. As in \cite{BaOk}, the author
established  one sufficient (a verification theorem) and one
necessary conditions of optimality.

In this paper, we aim at using convex analysis tools to prove
stochastic maximum principle for forward -backward stochastic
systems with jumps under partial information. More precisely, the
controlled systems are described by a  fully coupled nonlinear
multi-dimensional forward-backward stochastic differential equation
driven by a Poisson random measure and an independent Brownian
motion. Furthermore, all admissible control processes are required
to be adapted to a subfiltration of the filtration generated by the
underlying random measure and Brownian motion. The results obtained
in this paper can be considered as a generalization of \cite{Meng}
to the discontinuous case. Our paper covers the partial information
cases in \cite{BaOk} and \cite{Meng}. When there is no random
measure in the systems considered in our paper, the corresponding
stochastic optimal problem reduced to the case in \cite{Meng}. When
there is no backward system in the systems considered in our paper,
the corresponding stochastic optimal problem reduced to the case in
\cite{BaOk}.

  It is worth noting that in 2008,  Oksanal and Sulem \cite{OkSu}have investigated stochastic
  maximum principle for non-coupled one-dimensional forward-backward  differential equations with
  jumps. More precisely, the forward system does not couple with the
  backward system, only the backward system couples with the forward system.
  Compared with \cite{OkSu}, the system in our paper is fully coupled and
  multi-dimensional. So the system in \cite{OkSu} is a special case of the system in
  our paper.

The rest of this paper is organized as follows. In section 2, we
give our main assumptions and the statement of problem. In section
3,  we derive the main result, the sufficient maximum principle of
the stochastic optimal control problem under partial information.
Section 4 is devoted to the necessary optimality conditions.

Moreover, we refer to \cite{Wu1}\cite{Wu2} on the existence and
uniqueness of solutions to the fully coupled forward-backward
stochastic differential equations with jumps. There are already a
lot of literatures on the complete information maximum principle of
forward-backward stochastic differential systems. See
e.g.\cite{Xu95}\cite{ShWz06}\cite{ShWz07} and the references
therein.

\section{Statement of the optimal control problem}

 Let $(\Omega,{\cal F}, P)$ be a complete probability space and
$(\Omega,{\cal F}, \{{\cal F}_t\}_{t\geq 0}, P)$ be a filtered
probability space, where $\{{\cal F}_t\}_{t\geq 0}$  satisfies the
usual conditions, a right continuous increasing family of complete
sub $\sigma-$algebra of ${\cal F}$. let $\{{B_t}\}_{t\geq 0}$ be a
d-dimensional standard Brownian  motion in this space and let $\eta$
be a stationary ${\cal F}_t$-Poisson point process on a fixed
nonempty measurable subset $\mathscr{E}$ of $R^1$. We denote by
$\pi(de)$ the characteristic measure of $\eta$ and by $\tilde{N}(de,
dt)$ the counting measure induced by $\eta$. We assume that
$\pi(\mathscr{E})<\infty$,  we then define $N(de, dt):=\tilde{N}(de,
dt)-\pi(de)dt.$ We note that $N(de, dt)$ is a poisson martingale
measure with characteristic $\pi(de).$ We assume that $\{{\cal
F}_t\}_{t\geq 0}$ is the P-augmentation of the natural filtration
${\cal F}_t^{(W,~N)}$ defined by $\forall~ t\in (0,+\infty)$:
$$
{\cal F}_t^{(W,~N)}:=\sigma(W(s), 0\leq s\leq
t)\bigvee\sigma(\displaystyle N( A, (0, t]), 0\leq s\leq t,
 A\in \cal{B}(\mathscr{E})).$$

The following notation will be used in this paper:

$(\alpha,\beta):$ the inner product in Euclidean space $R^n, \forall
\alpha,\beta\in R^n.$

$|\alpha|=\sqrt{(\alpha,\alpha)}:$ the norm of Euclidean space
$R^n,\forall \alpha\in R^n.$

$(A,B)=tr(AB^T):$ the inner product in Euclidean space $R^{n\times
m},\forall A,B\in R^{n\times m}.$

$|A|=\sqrt{tr(AA^*)}:$ the norm of Euclidean space $R^{n\times
m},\forall A\in R^{n\times m}.   $

~~~~~~~~~Here we denote by A*,  the transpose of a matrix A

$(\alpha,\beta)_H:$ the inner product in Hilbert space $H, \forall
\alpha,\beta\in H.$

$|\alpha|_H=\sqrt{(\alpha,\alpha)_H}:$ the norm of Hilbert space
$H, \forall \alpha\in H.$

${\cal S}_{\cal F}^2(H):$ the Bananch space of $H$-valued ${\cal
F}_t$-adapted c\`{a}dl\`{a}g processes with the norm

~~~~~~~~~~~~$\sqrt{E\displaystyle \sup_{0\leq t\leq
T}|f(t)|_H^2}<~\infty.$

${\cal L}_{\cal F}^2(H):$ the Hilbert space of $H$-valued ${\cal
F}_t$-adapted processes with the norm

~~~~~~~~~~~~ $\sqrt{E\displaystyle\int_0^T|f(t)|_H^2dt}<~\infty.$

${\cal L}_\pi^2(H):$ the Hilbert space of H-valued measurable
  functions defined on the measure space

  ~~~~~~~~~~~ $(\mathcal {E}, \pi)$ with the norm
$\sqrt{\displaystyle\int_\mathcal {E}|r(z)|_H^2\pi(dz)}<~\infty.$

  $F_{\cal F}^{\pi, 2}(H):$ the Hilbert space of ${\cal L}_\pi^2(H)-$valued
${\cal F}_t$-predictable processes with the norm

~~~~~~~~~~~~~$\sqrt{E\displaystyle\int_0^T\int_\mathcal
{E}|r(t,e)|_H^2\pi (de)dt}<~\infty.$

$L^2(\Omega,{\cal F},P;H):$ the Hilbert space of $H$-valued norm
square integrable random variables on

~~~~~~~~~~~~~~~~~~$(\Omega,{\cal F},P).$

In this paper, we consider the controlled fully coupled nonlinear
forward-backward stochastic differential equations with jumps of the
form
\begin{equation}\label{eq:b1}
\displaystyle\left\{\begin{array}{lll}
dx_t&=&b(t,x_t,y_t,z_t,r_t(\cdot),v_t)dt+g(t,x_t,y_t,z_t,r_t(\cdot),v_{t})dB_t\\
&&+\displaystyle
\int_\mathscr{E}\sigma(t,x_{t-},y_{t-},z_t,r_t(\cdot),v_{t-},e){N}(de,
dt)\\
\displaystyle
dy_t&=&-f(t,x_t,y_t,z_t,r_t(\cdot),v_t)dt+z_tdB_t+\displaystyle
\int_\mathscr{E}r_t(e){N}(de, dt)\\\displaystyle
x(0)&=&a\\\displaystyle y(T)&=&\xi,

\end{array}
\right.
\end {equation}
where the coefficient are the maps as follows: $$
\begin{array}{ll}
\displaystyle &b(t, x, y, z, r(\cdot), v): [0,T]\times R^n\times
R^m\times
R^{m\times d}\times {\cal L}_\pi^2(R^m)\times {\cal U}\rightarrow R^n,\\
&g(t, x, y, z, r(\cdot), v): [0,T]\times R^n\times R^m\times
R^{m\times d}\times {\cal L}_\pi^2(R^m)\times {\cal U}\rightarrow
R^{n\times d},
\\
\displaystyle&\sigma(t, x, y, z, r(\cdot),v): [0,T]\times R^n\times
R^m\times R^{m\times d}\times {\cal L}_\pi^2(R^m)\times {\cal
U}\rightarrow R^{n},\\
\displaystyle &f(t,x, y, z,r(\cdot), v):[0,T]\times R^n\times
R^m\times R^{m\times d}\times {\cal L}_\pi^2(R^m)\times {\cal
U}\rightarrow R^m,\\
\end{array}
$$
Assume that b, g, $\sigma$, f are  Fr\`{e}chet differentiable with
respect to the variables  $(x,y,z,r(\cdot), v )$. For $\varphi=b, g,
\sigma, f$, we denote by $\nabla_x\varphi, \nabla_y\varphi,
\nabla_z\varphi, \nabla_u\varphi, \nabla_{r(\cdot)}\varphi$ the
Fr\`{e}chet derivatives with respect to x, y, z, u, $r(\cdot)$
respectively. $T> 0$ is a given constant, and $\xi$ is a given
random variable in $L^2(\Omega,{\cal F}_T, P; R^m)$. The process
$v_.=\{v_t(\omega), t\in[0, T], \omega\in \Omega\}$ in the system
\eqref{eq:b1}is the  control process, required to have values in a
given nonempty convex set ${\cal U}\subset R^k$ and required to be
c\`{a}dl\`{a}g and $\{\varepsilon_t \}_{t\geq 0}$ adapted, where
$\varepsilon_t\subseteq {\cal F}_t $ for all $t\in [0,~T]$ is a
given subfiltration representing the information available to the
controller at time t. For example, we could have
$\varepsilon_t={\cal F}_{(t-\delta)^+}$ for all  $t\in [0,~T],$
 where $\delta\geq 0$ is a fixed delay of information.

We shall define performance criterion by
\begin {equation}\label{eq:b2}
\displaystyle
J(v_.)=E\bigg[\displaystyle\int_0^Tl(t,x_t,y_t,z_t,r(\cdot),
v_t)dt+\phi(x_T)+h(y_0)\bigg],
\end {equation}
where
$$
\begin{array}{ll}
\displaystyle &l(t, x, y, z,r(\cdot), v): [0,T]\times R^n\times
R^m\times R^{m\times d}\times {\cal L}_\pi^2(R^m) \times {\cal
U}\rightarrow R,\\
\displaystyle&\phi( x): R^n \rightarrow R,\\
\displaystyle &h( y): R^m \rightarrow R,\\
\end{array}
$$
 are given Fr\`{e}chet  differential functions with respect to $(x,y,z, r(\cdot), v).$ We
 call the control process $v_.$  an admissible control if it gives rise
 to a a unique strong solution of the forward-backward stochastic
differential  equation \eqref{eq:b1} and the following condition
holds:
\begin{equation}\label{eq:b3}
 \displaystyle
E\bigg[\displaystyle\int_0^T|l(t,x_t,y_t,z_t,r(\cdot),v_t)|dt+|\phi(x_T)|+|h(y_0)|\bigg]<\infty.
\end{equation}
 The strong solution corresponding
to the admissible  control  $v_.$ is denoted by  $$\displaystyle
(x_.,y_.,z_., r_.(\cdot))=(x_.^{(v)},y_.^{(v)}, z_.^{(v)},
r_.^v(\cdot)) \in {\cal S}_{\cal F}^2(R^n)\times{\cal S}_{\cal
F}^2(R^m)\times {\cal L}_{\cal F}^2(R^{m\times d}) \times F_{\cal
F}^{\pi, 2}(R^m).$$
 The set of all admissible
controls is denoted by ${\cal A}$. If $v_.\in {\cal A }$ and
$(x_.,y_.,z_., r_.(\cdot))=(x_.^{(v)},y_.^{(v)}, z_.^{(v)},
r_.^v(\cdot))$ is the corresponding strong solution of
\eqref{eq:b1}, we call $(v_.; x_.,y_.,z_., r_.(\cdot))$ an
admissible pair.

The partial information optimal control problem amounts to
determining an admissible control  $u_.\in {\cal A }$ such that
\begin{equation}\label{eq:b4}
J(u_.)=\displaystyle\inf_{v_.\in {\cal A}}J(v_.).
\end{equation} Such controls $u_.$ are called optimal control. If
$(x_.,y_.,z_., r(\cdot))=(x_.^{(u)},y_.^{(u)}, z_.^{(u)},
r_.^{(u)}(\cdot))$ is the corresponding strong solution of
\eqref{eq:b1}, then $(u_.;  x_., y_., z_., r_.(\cdot))$ is called an
optimal pair.
\section{A Partial Information Sufficient Maximum Principle}
In this section we want to study the  sufficient maximum principle
for the partial information optimal control problem
\eqref{eq:b1}\eqref{eq:b2}\eqref{eq:b3}\eqref{eq:b4}. To this end,
we need the following two lemmas. The main results of this paper are
the following

\begin{lem}
\textbf{(Integration by Parts)} Suppose that the processes
$Y^{(1)}(t)$ and $Y^{(2)}(t)$ are given by
$$
dY^{(j)}(t)=b^{(j)}(t)dt+g^{(j)}(t)dB_t+\displaystyle
\int_\mathscr{E} \sigma^{(j)}(t, e){N}(de, dt),
Y^{(j)}(0)=y^{(j)}\in R^n,  j=1.2,$$ where $b^{(j)}\in {\cal
L}_{F}^2(R^n), g^{(j)}\in {\cal L }_{\cal F}^2(R^{n\times d}),
\sigma^{(j)}\in F^{\pi, 2}_{\cal F}(R^n).$ Then,
$$\begin {array} {ll}
E[Y^{(1)}(T)Y^{(2)}(T)]=&y_1y_2+E\bigg[\displaystyle
\int^T_0\bigg(Y^{(1)}(t-), dY^{(2) }(t)\bigg)\bigg]
+E\bigg[\displaystyle \int^T_0\bigg(dY^{(1)}(t), Y^{(2)
}(t-)\bigg)\bigg]\\& +E\bigg[\displaystyle \int^T_0\bigg(g^{(1)}(t),
g^{(2)}(t)\bigg)dt\bigg] +E\bigg[\displaystyle \int^T_0\displaystyle
\int_\mathscr{E}\bigg(\sigma^{(1)}(t,e),
\sigma^{(2)}(t,e)\bigg)\pi(de)dt\bigg].
\end {array}$$
\end{lem}
The proof can be obtained directly from the It$\grave{o}$ formula
with jumps (Theorem 1.16 in \cite{OkSu})
\begin{lem}
  Let f be a Fr\`{e}chet differentiable on Hilbert space H. And let C be a convex subset of H.
       Then f is convex on C if and only if
       $$f(x)-f(x_0)\geq \bigg(\nabla_{x_0} f,~~x-x_0\bigg)_H$$ for all $x, x_0 \in C.$

 The proof can be obtained directly from Theorem 4.1.1. in \cite{BaLe}.
\end{lem}
In our setting  the Hamiltonian function $H:[0,T]\times R^n\times
R^m\times R^{m\times d}\times {\cal L}_\pi^2(R^m)\times {\cal
U}\times R^n\times R^{n\times d}\times{\cal L}_\pi^2(R^n)\times
R^m\rightarrow R$ gets the following form:
\begin {equation}
\begin{array}{ll}
\displaystyle H(t,x,y,z,r(\cdot),v,p,q,\beta(\cdot),k)=&\bigg(k,
-f(t,x,y,z,r(\cdot),v)\bigg)+\bigg(p,
b(t,x,y,z,r(\cdot),v)\bigg)\\\displaystyle
&+\bigg(q,\sigma(t,x,y,z,r(\cdot), v)\bigg)
+l(t,x,y,z,r(\cdot),v)\\&+\displaystyle\int_\mathscr{E}\bigg(\beta(\bar{e}),
\sigma(t,x,y,z,r(\cdot),v,\bar{e})\bigg)\pi(d\bar{e}) .
\end{array}
\end {equation}The  adjoint  equation which fits into  the system (2.1) corresponding
   to the given admissible pair $(v_.;  x_., y_., z_., r_.(\cdot))$ is given by the following
   forward-backward stochastic differential equation:
   \begin {equation}\label{eq:c2}
  \left\{\begin{array}{lll}
  \displaystyle
  dk_t&=&-\nabla_yH(t,x_t,y_t,z_t,r_t(\cdot),v_t,p_t,q_t,\beta(\cdot),k_t)dt\\
  &&-\nabla_zH(t,x_t,y_t,z_t,r_t(\cdot),v_t,p_t,q_t,\beta(\cdot),k_t)dB_t,\\
  &&-\displaystyle\int_\mathscr{E}\nabla_{r(\cdot)}H(t,x_t,y_t,z_t,r(\cdot), v_t,p_t,q_t,\beta(\cdot),k_t)N(de,dt)\\\displaystyle
dp_t&=&-\nabla_x H(t,x_t,y_t,z_t,r_t(\cdot),v_t,p_t,q_t,\beta(\cdot),k_t)dt+q_tdB_t+\displaystyle\int_\mathscr{E}\beta(e)N(de,dt),\\
 k_0&=&-\nabla_yh(y_0),
~~~p_T=\nabla _x\varphi(x_T),
  \end{array}
  \right.
  \end {equation}
where $(p(t),q(t),\beta(t,\cdot),k(t))\in R^n\times R^{n\times
d}\times{\cal L}_\pi^2(R^n)\times R^m$ are the unknown processes.

We now coming to a verification theorem for the optimal control
problem \eqref{eq:b1}\eqref{eq:b2}\eqref{eq:b3}\eqref{eq:b4}.

    \begin{thm}\textbf{(The Sufficient Maximum Principle)}\label{thm:c1}

Let $(\hat{u_.}; \hat{x_.},\hat{y_.},\hat{z_.}, \hat{r}_.(\cdot))$
be an admissible pair and suppose that there exists a strong
solution $(\hat{p}_.,\hat{q}_.,\hat{\beta}_.(\cdot), \hat{k}_.)$ of
the corresponding adjoint forward-backward equation \eqref{eq:c2}.
 Assume that the following conditions are satisfied. For arbitrary admissible
control $({v_.}; {x_.^{(v)}},{y_.^{(v)}},{z_.^{(v)}},
r_.^{(v)}(\cdot))$,  we have
\begin{equation}
 \displaystyle E\displaystyle\int_0^T(\hat{x}_t-x_t^{(v)})^*\hat{q}_t\hat{q}_t^*(\hat{x}_t
 -x_t^{(v)})dt<+\infty,
\end {equation}
\begin{equation}
 \displaystyle E\displaystyle\int_0^T(\hat{x}_t-x_t^{(v)})^*\int_\mathscr{E}
 \hat{\beta}_t(e)\hat{\beta}^*_t(e)\pi(de)(\hat{x}_t-x_t^{(v)})dt<+\infty,
\end {equation}
\begin {equation}
  E\displaystyle\int_0^T(\hat{y}_t-y_t^{(v)})^*( \nabla_zH \nabla_zH^*)(t,\hat {x}_t,
  \hat{y}_t, \hat{z}_t, \hat{r}(\cdot),\hat{v}_t,\hat{p}_t,
 \hat{q}_t,\hat{\beta}(\cdot),\hat{k}_t)(\hat{y}_t-y_t^{(v)})dt<+\infty,
\end {equation}
\begin {equation}
  E\displaystyle\int_0^T(\hat{y}_t-y_t^{(v)})^*\int_\mathscr{E}( \nabla_{r(\cdot)}H
   \nabla_{r(\cdot)}H^*)(t,\hat {x}_t,\hat{y}_t, \hat{z}_t, \hat{r}(\cdot),\hat{v}_t,\hat{p}_t,
 \hat{q}_t,\hat\beta(\cdot),\hat{k}_t)\pi(dz)(\hat{y}_t-y_t^{(v)})dt<+\infty,
\end {equation}
\begin {equation}
E\displaystyle\int_0^T\hat{p}_t^*(gg^*)(t,x_t^{(v)},y_t^{(v)},r_t^{(v)}(\cdot),
z_t^{(v)},v_t)\hat{p}_tdt<+\infty,
\end {equation}
\begin {equation}
E\displaystyle\int_0^T\hat{p}_t^*\int_\mathscr{E}(\sigma\sigma^*)(t,x_t^{(v)},y_t^{(v)},z_t^{(v)},r_t^{(v)}
(\cdot),v_t, e)\pi(d e)\hat{p}_tdt<+\infty,
\end {equation}
\begin {equation}
  E\displaystyle\displaystyle\int_0^T\hat{k}_t^*(z_t^{(v)}z_t^{(v)~*})
  \hat {k}_tdt<+\infty,
\end {equation}
\begin {equation}
  E\displaystyle\displaystyle\int_0^T\hat{k}_t^*\int_\mathscr E(r_t^{(v)}
  (e)r_t^{(v)~*}(e))\pi(de)\hat {k}_tdt<+\infty,
\end {equation}
\begin {equation}
  E\displaystyle\int_0^T|H_u(t,\hat{x}_t,\hat{y}_t,
  \hat{z}_t,\hat{r}_{\cdot},\hat{u}_t, \hat{p}_t,\hat{q}_t,\hat{\beta}{(\cdot)},\hat{k}_t)|^2<+\infty.
\end {equation}
ensuring that the integrales with repect to the Brorwnian motion B
and the compensated jump parts indeed have zero mean.
 Moreover, suppose that for all $t\in
[0,T]$, $H(t,x,y,z,r_t(\cdot),
v,\hat{p}_t,\hat{q}_t,\hat{\beta}_t(\cdot), \hat{k}_t)$ is convex in
$(x,y,z,r(\cdot), v)$ , and $h(y)$ is convex in $y$ and $\phi(x)$ is
convex in $x$, and the following partial information maximum
condtion holds
\begin {equation}
E\bigg[H(t,\hat{x}_t,\hat{y}_t,\hat{z}_t,\hat{r}(\cdot),\hat{u}_t,\hat{p}_t,
\hat{q}_t,\hat{\beta}(\cdot),\hat{k}_t)\bigg|\varepsilon_t\bigg]=
\displaystyle min_{v\in
U}E\bigg[H(t,\hat{x}_t,\hat{y}_t,\hat{z}_t,\hat{r}(\cdot),v,\hat{p}_t,
\hat{q}_t,\hat{\beta}(\cdot),\hat{k}_t)\bigg|\varepsilon_t\bigg].
\end {equation} Then $\hat{u}_.$ is a partial information optimal
control.
\end{thm}

\begin{proof}
  Let $(v_.; x_.,y_.,z_.,r(\cdot))=(v_.;
x_.^{(v)},y_.^{(v)},r_.^{(v)}(\cdot),z_.^{(v)})$ be an arbitrary
admissible pair. It follows from the definition of the performance
functional \eqref{eq:b2} that
\begin{equation}
\begin{array}{ll}
 J(v(\cdot))-J(\hat{u}(\cdot))&=
E\displaystyle\int_0^T\bigg[l(t,x_t,y_t,z_t,r_t(\cdot),v_t)-l(t,\hat{x}_t,\hat{y}_t,\hat{z}_t,
\hat{r}_t{(\cdot)},\hat{u}_t)\bigg]dt
\\&~~+E\bigg[\phi(x_T)-\phi(\hat{x}_T)\displaystyle\bigg]
+E\bigg[h(y_0)-h(\hat{y}_0)\bigg]\\&=I_1+I_2,
\end{array}
\end{equation}
where \begin{equation} \displaystyle
I_1=E\int_0^T\bigg[l(t,x_t,y_t,z_t,r_t(\cdot),v_t)-l(t,\hat{x}_t,\hat{y}_t,\hat{z}_t,\hat{r}_t(\cdot),
 \hat{u}_t)\bigg]dt,
\end{equation}
\begin{equation}
\begin{array}{ll}
\displaystyle
I_2=E\bigg[\phi(x_T)-\phi(\hat{x}_T)\bigg]+E[h(y_0)-h(\hat{y}_0)].
\end{array}
\end{equation}
 Using Convexity of $\varphi$ and h, Lemma 3.2 and Lemma 3.1, we have
$$
\begin{array}{lll}
I_2&=&E\bigg[\phi(x_T)-\phi(\hat{x}_T)\bigg]+E\bigg[h(y_0)-h(\hat{y}_0)\bigg]\\
&\geq&
E\bigg[\bigg(\nabla_x\phi(\hat{x}_T),x_T-\hat{x}_T\bigg)\bigg]+E\bigg[\bigg(\nabla_yh(\hat{y}_0),
y_0-\hat{y}_0\bigg)\bigg]
\\\displaystyle
&=&E\bigg[\bigg(\hat{p}_T, x_T-\hat{x}_T\bigg)-\bigg(\hat{p}_0,x_0
-\hat{x}_0\bigg)+\bigg(\hat{k}_T,
y_T-\hat{y}_T\bigg)-\bigg(\hat{k}_0, y_0-\hat{y}_0\bigg)\bigg]\\
\displaystyle
&=&E\bigg[\displaystyle\int_0^T\bigg(x_{t-}-\hat{x}_{t-}, d\hat{p}_t
\bigg)\bigg]+E\bigg[\int_0^T
\bigg(\hat{p}_{t-}, d(x_t-\hat{x}_t)\bigg)\bigg]\\
\displaystyle &&+ E\bigg[\displaystyle\int_0^T\bigg(\hat{q}_t,
g(t,x_t,y_t,z_t,r_t(\cdot),v_t)-g(t,\hat{x}_t,\hat{y}_t,\hat{z}_t,\hat{r}_t(\cdot),\hat{u}_t)\bigg)\bigg]dt
\\
\displaystyle &&+
E\bigg[\displaystyle\int_0^T\displaystyle\int_\mathscr{E}\bigg(\hat{\beta}_t(e),
\sigma(t,x_t,y_t,z_t,r_t(\cdot),v_t,e)-\sigma(t,\hat{x}_t,\hat{y}_t,\hat{z}_t,\hat{r}_t(\cdot),
\hat{u}_t,e)\bigg)\bigg]\pi(de)dt
\\
\end{array}
$$
$$
\begin{array}{lll}
&&+E\bigg[\displaystyle\int_0^T\bigg(y_{t-}-\hat{y}_{t-},
d\hat{k}_t\bigg)\bigg] \displaystyle
+E\bigg[\displaystyle\int_0^T\bigg(\hat{k}_{t-},d(y_t-\hat{y}_t)\bigg)\bigg]
\\&&+E\bigg[\displaystyle\int_0^T\bigg(-\nabla_zH(t,\hat{x}_t,\hat{y}_t,\hat{z}_t,\hat{r}(\cdot),\hat{u}_t,\hat{p}_t,
\hat{q}_t,\hat{\beta}(\cdot),\hat{k}_t),z_t-\hat{z}_t\bigg)dt
\\&&+E\bigg[\displaystyle\int_0^T\displaystyle\int_\mathscr{E}\bigg(-\nabla_{r(\cdot)}H(t,\hat{x}_t,
\hat{y}_t,\hat{z}_t,\hat{r}(\cdot),\hat{u}_t,\hat{p}_t,
\hat{q}_t,\hat{\beta}(\cdot),\hat{k}_t),r_t(e)-\hat{r}_t(e)\bigg)\pi(de)dt\bigg]
\\
&=&E\bigg[\displaystyle\int_0^T\displaystyle\bigg(-\nabla_{x}H(t,\hat{x}_t,\hat{y}_t,\hat{z}_t,\hat{r}(\cdot),
\hat{u}_t,\hat{p}_t,
\hat{q}_t,\hat{\beta}(\cdot),\hat{k}_t),x_t-\hat{x}_t\bigg)dt\bigg]
\displaystyle
\\&&+E\bigg[\displaystyle\int_0^T\displaystyle\bigg(-\nabla_{y}H(t,\hat{x}_t,\hat{y}_t,\hat{z}_t,\hat{r}
(\cdot),\hat{u}_t,\hat{p}_t,
\hat{q}_t,\hat{\beta}(\cdot),\hat{k}_t),y_t-\hat{y}_t\bigg)dt\bigg]
\displaystyle
\\
&&+E\bigg[\displaystyle\int_0^T\displaystyle\bigg(-\nabla_{z}H(t,\hat{x}_t,\hat{y}_t,\hat{z}_t,\hat{r}
(\cdot),\hat{u}_t,\hat{p}_t,
\hat{q}_t,\hat{\beta}(\cdot),\hat{k}_t),z_t-\hat{z}_t\bigg)dt\bigg]
\displaystyle
\\
&&+E\bigg[\displaystyle\int_0^T\displaystyle\int_\mathscr{E}\bigg(-\nabla_{r(\cdot)}H(t,\hat{x}_t,
\hat{y}_t,\hat{z}_t,\hat{r}(\cdot),\hat{u}_t,\hat{p}_t,
\hat{q}_t,\hat{\beta}(\cdot),\hat{k}_t),r_t(e)-\hat{r}_t(e)\bigg)\pi(de)dt\bigg]
\displaystyle
\\
&&+E\bigg[\displaystyle\int_0^T\bigg(\hat{p}_t,b(t,x_t,y_t,z_t,r_t(\cdot),v_t)-b(t,\hat{x}_t,\hat{y}_t,\hat{z}_t,
\hat{r}_t(\cdot),\hat{u}_t)\bigg)dt\\
&&+E\bigg[\displaystyle\int_0^T\bigg(\hat{q}_t,g(t,x_t,y_t,z_t,r_t(\cdot),v_t)-g(t,\hat{x}_t,\hat{y}_t,\hat{z}_t,
\hat{r}_t(\cdot),\hat{u}_t)\bigg)dt
\\&&+E\bigg[\displaystyle\int_0^T\displaystyle\int_\mathscr{E}\bigg(\hat{\beta}_t(e),
\sigma(t,x_t,y_t,z_t,r_t(\cdot),v_t,e)-\sigma(t,\hat{x}_t,\hat{y}_t,\hat{z}_t,
\hat{r}_t(\cdot),\hat{u}_t,e)\bigg)\pi(de)dt
\\&&+E\bigg[\displaystyle\int_0^T\bigg(\hat{k}_t,-(f(t,x_t,y_t,z_t,r_t(\cdot),v_t)-f(t,\hat{x}_t,\hat{y}_t,\hat{z}_t,
\hat{r}_t(\cdot),\hat{u}_t))\bigg)dt
\\
&=&J_1+J_2,
\end{array}
$$
where
$$
\begin{array}{ll}
J_1=&E\bigg[\displaystyle\int_0^T\displaystyle\bigg(-\nabla_{x}H(t,\hat{x}_t,\hat{y}_t,\hat{z}_t,
\hat{r}(\cdot),\hat{u}_t,\hat{p}_t,
\hat{q}_t,\hat{\beta}(\cdot),\hat{k}_t),x_t-\hat{x}_t\bigg)dt\bigg]
\displaystyle
\\&+E\bigg[\displaystyle\int_0^T\displaystyle\bigg(-\nabla_{y}H(t,\hat{x}_t,\hat{y}_t,\hat{z}_t,
\hat{r}(\cdot),\hat{u}_t,\hat{p}_t,
\hat{q}_t,\hat{\beta}(\cdot),\hat{k}_t),y_t-\hat{y}_t\bigg)dt\bigg]
\displaystyle
\\
&+E\bigg[\displaystyle\int_0^T\displaystyle\bigg(-\nabla_{z}H(t,\hat{x}_t,\hat{y}_t,\hat{z}_t,
\hat{r}(\cdot),\hat{u}_t,\hat{p}_t,
\hat{q}_t,\hat{\beta}(\cdot),\hat{k}_t),z_t-\hat{z}_t\bigg)dt\bigg]
\displaystyle
\\
&+E\bigg[\displaystyle\int_0^T\displaystyle\int_\mathscr{E}\bigg(-\nabla_{r(\cdot)}H(t,\hat{x}_t,
\hat{y}_t,\hat{z}_t,\hat{r}(\cdot),\hat{u}_t,\hat{p}_t,
\hat{q}_t,\hat{\beta}(\cdot),\hat{k}_t),r_t(e)-\hat{r}_t(e)\bigg)\pi(de)dt\bigg]
\displaystyle,
\end{array}
$$
$$
\begin{array}{ll}
J_2=&E\bigg[\displaystyle\int_0^T\bigg(\hat{p}_t,b(t,x_t,y_t,z_t,r_t(\cdot),v_t)-b(t,\hat{x}_t,\hat{y}_t,\hat{z}_t,
\hat{r}_t(\cdot),\hat{u}_t)\bigg)dt\bigg]\\
&+E\bigg[\displaystyle\int_0^T\bigg(\hat{q}_t,g(t,x_t,y_t,z_t,r_t(\cdot),v_t)-g(t,\hat{x}_t,\hat{y}_t,\hat{z}_t,
\hat{r}_t(\cdot),\hat{u}_t)\bigg)dt\bigg]
\\&+E\bigg[\displaystyle\int_0^T\displaystyle\int_\mathscr{E}\bigg(\hat{\beta}_t(e),
 \sigma(t,x_t,y_t,z_t,r_t(\cdot),v_t,e)-\sigma(t,\hat{x}_t,\hat{y}_t,\hat{z}_t,
\hat{r}_t(\cdot),\hat{u}_t,e)\bigg)\pi(de)dt\bigg]
\\&+E\bigg[\displaystyle\int_0^T\bigg(\hat{k}_t,-(f(t,x_t,y_t,z_t,r_t(\cdot),v_t)-f(t,\hat{x}_t,\hat{y}_t,\hat{z}_t,
\hat{r}_t(\cdot),\hat{u}_t))\bigg)dt\bigg],
\end{array}
$$
and we have used the assumptions (3.3)-(3.11)which ensure that the
stochastic integrals with respect to the Brownian motion and the
Poisson random measure have zero expectation.

By the definition of the Hamiltonian function $H$ (noting (3.1)) and
the definition of $I_1$ (noting (3.14)), we have
\begin{equation}
\begin{array}{ll}
I_1=&E\displaystyle\int_0^T\bigg[l(t,x_t,y_t,z_t,r_t(\cdot),v_t)-l(t,\hat{x}_t,
\hat{y}_t,\hat{z}_t,\hat{r}_t(\cdot),\hat{u}_t)\bigg]dt\\
~~~=&E\displaystyle\int_0^T\bigg[H(t,x_t,y_t,z_t,r_t(\cdot),v_t,\hat{p}_t,\hat{q}_t,\hat{\beta}_t(\cdot),\hat{k}_t)
-H(t,\hat{x}_t,\hat{y}_t,\hat{z}_t,\hat{r}_t(\cdot),\hat{u}_t,\hat{p}_t,\hat{q}_t,\hat{\beta}(\cdot),\hat{k}_t)\bigg]dt
\\&-E\bigg[\displaystyle\int_0^T\bigg(\hat{p}_t,b(t,x_t,y_t,z_t,r_t(\cdot),v_t)-b(t,\hat{x}_t,\hat{y}_t,\hat{z}_t,
\hat{r}_t(\cdot),\hat{u}_t)\bigg)dt\bigg]\\
&-E\bigg[\displaystyle\int_0^T\bigg(\hat{q}_t,g(t,x_t,y_t,z_t,r_t(\cdot),v_t)-g(t,\hat{x}_t,\hat{y}_t,\hat{z}_t,
\hat{r}_t(\cdot),\hat{u}_t)\bigg)dt\bigg]
\\&-E\bigg[\displaystyle\int_0^T\displaystyle\int_\mathscr{E}\bigg(\hat{\beta}_t(e),
\sigma(t,x_t,y_t,z_t,r_t(\cdot),v_t,e)-\sigma(t,\hat{x}_t,\hat{y}_t,\hat{z}_t,
\hat{r}_t(\cdot),\hat{u}_t,e)\bigg)\pi(de)dt\bigg]
\\&-E\bigg[\displaystyle\int_0^T\bigg(\hat{k}_t,-(f(t,x_t,y_t,z_t,r_t(\cdot),v_t)-f(t,\hat{x}_t,\hat{y}_t,\hat{z}_t,
\hat{r}_t(\cdot),\hat{u}_t))\bigg)dt\bigg],
\\
~~~=&J_3-J_2,
\end{array}
\end{equation}
where
\begin{equation}
J_3=E\displaystyle\int_0^T\bigg[H(t,x_t,y_t,z_t,r_t(\cdot),v_t,\hat{p}_t,\hat{q}_t,\hat{\beta}_t(\cdot),\hat{k}_t)
-H(t,\hat{x}_t,\hat{y}_t,\hat{z}_t,\hat{r}_t(\cdot),
\hat{u}_t,\hat{p}_t,\hat{q}_t,\hat{\beta}(\cdot),\hat{k}_t)\bigg]dt.
\end{equation}
From the convexity of
$H(t,x,y,z,r(\cdot),v,\hat{p}_t,\hat{q}_t,\hat{\beta}_t(\cdot),\hat{k}_t))$
with respect to $(x,y,z,r(\cdot),v )$ and Lemma 3.2, we have
\begin{equation}
\begin{array}{ll}
&H(t,x_t,y_t,z_t,r_t(\cdot),v_t,\hat{p}_t,\hat{q}_t,\hat{\beta}_t(\cdot),\hat{k}_t)
-H(t,\hat{x}_t,\hat{y}_t,\hat{z}_t,\hat{r}_t(\cdot),\hat{u}_t,\hat{p}_t,\hat{q}_t,\hat{\beta}(\cdot),\hat{k}_t)
\\ \geq&
\bigg(\nabla_{x}H(t,\hat{x}_t,\hat{y}_t,\hat{z}_t,\hat{r}(\cdot),\hat{u}_t,\hat{p}_t,
\hat{q}_t,\hat{\beta}(\cdot),\hat{k}_t),x_t-\hat{x}_t\bigg)
\displaystyle
\\&+\bigg(\nabla_{y}H(t,\hat{x}_t,\hat{y}_t,\hat{z}_t,\hat{r}(\cdot),\hat{u}_t,\hat{p}_t,
\hat{q}_t,\hat{\beta}(\cdot),\hat{k}_t),y_t-\hat{y}_t\bigg)
\displaystyle
\\
&+\bigg(\nabla_{z}H(t,\hat{x}_t,\hat{y}_t,\hat{z}_t,\hat{r}(\cdot),\hat{u}_t,\hat{p}_t,
\hat{q}_t,\hat{\beta}(\cdot),\hat{k}_t),z_t-\hat{z}_t\bigg)
\displaystyle
\\
&+\displaystyle\int_\mathscr{E}\bigg(\nabla_{r(\cdot)}H(t,\hat{x}_t,\hat{y}_t,\hat{z}_t,
\hat{r}(\cdot),\hat{u}_t,\hat{p}_t,
\hat{q}_t,\hat{\beta}(\cdot),\hat{k}_t),r_t(e)-\hat{r}_t(e)\bigg)\pi(de)
\displaystyle
\\&+\bigg(\nabla_{v}H(t,\hat{x}_t,\hat{y}_t,\hat{z}_t,\hat{r}(\cdot),\hat{u}_t,\hat{p}_t,
\hat{q}_t,\hat{\beta}(\cdot),\hat{k}_t),v_t-\hat{u}_t\bigg).
\end{array}
\end{equation}
Since $v\rightarrow
E[H(t,\hat{x}_t,\hat{y}_t,\hat{z}_t,\hat{r}(\cdot),v,\hat{p}_t,
\hat{q}_t,\hat{\beta}(\cdot),\hat{k}_t)|\varepsilon_t],~v\in{\cal
U}$ is minimal for $\hat{u}_t$ and  $v_t,\hat{u}_t$ are
$\varepsilon_t$-measurable, we can get by (3.8)
\begin{equation}
\begin{array}{ll}
&
\bigg(E\bigg[\nabla_{v}H(t,\hat{x}_t,\hat{y}_t,\hat{z}_t,\hat{r}(\cdot),\hat{u}_t,\hat{p}_t,
\hat{q}_t,\hat{\beta}(\cdot),\hat{k}_t)\bigg|
\varepsilon_t\bigg], v_t-\hat{u}_t\bigg)\\
=&E\bigg[\bigg(\nabla_{v}H(t,\hat{x}_t,\hat{y}_t,\hat{z}_t,\hat{r}(\cdot),\hat{u}_t,\hat{p}_t,
\hat{q}_t,\hat{\beta}(\cdot),\hat{k}_t),v_t-\hat{u}_t\bigg)
\bigg|\varepsilon_t\bigg]\geq 0.
\end{array}
\end{equation}
Hence combining (3.17), (3.18) and (3.19), we obtain
\begin{equation}
\begin{array}{lll}
J_3&\geq &
E\bigg[\displaystyle\int_0^T\displaystyle\bigg(\nabla_{x}H(t,\hat{x}_t,\hat{y}_t,\hat{z}_t,
\hat{r}(\cdot),\hat{u}_t,\hat{p}_t,
\hat{q}_t,\hat{\beta}(\cdot),\hat{k}_t),x_t-\hat{x}_t\bigg)dt\bigg]
\displaystyle
\\&&+E\bigg[\displaystyle\int_0^T\displaystyle\bigg(\nabla_{y}H(t,\hat{x}_t,\hat{y}_t,\hat{z}_t,
\hat{r}(\cdot),\hat{u}_t,\hat{p}_t,
\hat{q}_t,\hat{\beta}(\cdot),\hat{k}_t),y_t-\hat{y}_t\bigg)dt\bigg]
\displaystyle
\\
&&+E\bigg[\displaystyle\int_0^T\displaystyle\bigg(\nabla_{z}H(t,\hat{x}_t,\hat{y}_t,
\hat{z}_t,\hat{r}(\cdot),\hat{u}_t,\hat{p}_t,
\hat{q}_t,\hat{\beta}(\cdot),\hat{k}_t),z_t-\hat{z}_t\bigg)dt\bigg]\\
&&+E\bigg[\displaystyle\int_0^T\displaystyle\int_\mathscr{E}\bigg(\nabla_{r(\cdot)}H(t,\hat{x}_t,\hat{y}_t,\hat{z}_t,
\hat{r}(\cdot),\hat{u}_t,\hat{p}_t,
\hat{q}_t,\hat{\beta}(\cdot),\hat{k}_t),r_t(e)-\hat{r}_t(e)\bigg)\pi(de)dt\bigg],\\
&=&J_1.
\end{array}
\end{equation}
Therefore, it follows from (3.10), (3.15), (3.16) and (3.20) that
$$
J(v_.)-J(\hat{u}_.)=I_1+I_2=(J_3-J_2)+I_2\geq
(J_1-J_2)+(-J_1+J_2)=0.
$$
Since $v_.\in{\cal A}$ is arbitrary, we conclude that $\hat{u}_.$ is
optimal control. The proof of Theorem 3.1 is completed.
\end{proof}
\section {A Partial Information Necessary Maximum Principle }
In this section, we give a necessary maximum principle for the
partial information  stochastic optimal control problem
\eqref{eq:b1}\eqref{eq:b2}\eqref{eq:b3}\eqref{eq:b4}. To this end,
we adopt a similar strategy as in \cite{BaOk} and \cite{Meng}.

In addition to the assumption in Section 2, we shall now assume the
following:

$(H_1)$ For all $t, r$ such that $0\leq t<t+r\leq T,$ and all
bounded $\varepsilon_t$-measurable $\alpha=\alpha(\omega),$ the
control $\theta_s:=(0,\cdots, \theta^i_s,0,\cdots,0)\in {\cal
U}\subset R^k,$with $\theta^i_s=:\alpha_i\chi_{[t,~t+r]}(s),~~s\in
[0,T]$ belongs to ${\cal A},  i=1,2,\cdots, k$.

$(H_2)$ For given $u_., \theta_.\in {\cal A}$,  with $\theta_.$
bounded, there exists $\delta>0$ such that $u_.+y\theta_.\in {\cal
A},$ for all $y\in (-\delta,\delta)$.

For given $u_., \theta_.\in {\cal A}$ with $\theta_.$ bounded, we
define the processes $(X_t^{1},Y_t^{1},Z_t^{1}, R_t^{1} )_{t\geq 0}$
by
\begin{equation}
\begin{array}{ll}
\displaystyle X_t^{1}&:=X_t^{(u, ~ \theta)}:=\displaystyle\frac{d}{dy}x^{(u+y\theta)}_t|_{y=0},\\
\displaystyle Y_t^{1}&:=Y_t^{(u,  ~\theta)}:=\displaystyle\frac{d}{dy}y^{(u+y\theta)}_t|_{y=0},\\
\displaystyle Z_t^{1}&:=Z_t^{(u,
~\theta)}:=\displaystyle\frac{d}{dy}z^{(u+y\theta)}_t|_{y=0},\\
\displaystyle R_t^{1}&:=R_t^{(u,
~\theta)}:=\displaystyle\frac{d}{dy}r^{(u+y\beta)}_t|_{y=0}.
\end{array}
\end{equation}
Note  that $(X_t^{1},Y_t^{1},Z_t^{1},R_t^{1})_{t\geq 0}$ satisfies
the following linear forward-backward stochastic differential
equation
\begin{equation}
\displaystyle\left\{\begin{array}{ll}
dX_t^{1}=&\bigg[\nabla_xb(t,x_t,y_t,z_t,r_t(\cdot),u_t)X_t^{1}+\nabla_yb(t,x_t,y_t,z_t,r_t(\cdot),u_t)Y_t^{1}
\\\displaystyle
&+\nabla_zb(t,x_t,y_t,z_t,r_t(\cdot),u_t)Z_t^{1}
+\displaystyle\int_\mathscr{E}\nabla_{r(\cdot)}b(t,x_t,y_t,z_t,r_t(\cdot),u_t)R_t(\bar{e})\pi(d\bar{e})
\\&+\nabla_vb(t,x_t,y_t,z_t,r_t(\cdot),u_t)\theta_t\bigg]dt\displaystyle
+\bigg[\nabla_xg(t,x_t,y_t,z_t,r_t(\cdot),u_t)X_t^{1}\\
&+\nabla_yg(t,x_t,y_t,z_t,r_t(\cdot),u_t)Y_t^{1} \displaystyle
+\nabla_zg(t,x_t,y_t,z_t,r_t(\cdot),u_t)Z_t^{1}\\
&+\displaystyle\int_\mathscr{E}\nabla_{r(\cdot)}g(t,x_t,y_t,z_t,r_t(\cdot),u_t)R_t(\bar{e})\pi(d\bar{e})
+\nabla_vg(t,x_t,y_t,z_t,r_t(\cdot),u_t)\theta_t\bigg]dB_t
\\\displaystyle
&+\displaystyle\int_\mathscr{E}\bigg[\nabla_x\sigma(t,x_t,y_t,z_t,r_t(\cdot),u_t,e)X_t^{1}
+\nabla_y\sigma(t,x_t,y_t,z_t,r_t(\cdot),u_t,e)Y_t^{1}
\\\displaystyle
&+\nabla_z\sigma(t,x_t,y_t,z_t,r_t(\cdot),u_t, e)Z_t^{1}
+\displaystyle\int_\mathscr{E}\nabla_{r(\cdot)}\sigma(t,x_t,y_t,z_t,r_t(\cdot),u_t,e)R_t(\bar{e})\pi(d\bar{e})
\\&+\nabla_v\sigma(t,x_t,y_t,z_t,r_t(\cdot),u_t,e)\theta_t\bigg]N(de,dt)
\\dY_t^{1}=&-\bigg[\nabla_xf(t,x_t,y_t,z_t,r_t(\cdot),u_t)X_t^{1}+\nabla_yf(t,x_t,y_t,z_t,r_t(\cdot),u_t)Y_t^{1}\\
&+\nabla_zf(t,x_t,y_t,z_t,r_t(\cdot),u_t)Z_t^{1}
+\displaystyle\int_\mathscr{E}\nabla_{r(\cdot)}f(t,x_t,y_t,z_t,r_t(\cdot),u_t)R_t^{1}(\bar{e})\pi(d\bar{e})
\\&+\nabla_vf(t,x_t,y_t,z_t,r_t(\cdot),u_t)\theta_t\bigg]dt+Z_t^{1}dB_t
+\displaystyle\int_\mathscr{E} R_t^{1}(e) N(de,dt)\\
\displaystyle X_0^{1}=0,&Y_T^{1}=0,
\end{array}
\right.
\end{equation}

where $({x}_t, {y}_t, {z}_t,
r_t(\cdot))=(x_t^{({u})},{y}_t^{({u})},{z}_t^{({u})},{r}_t^{({u})}(\cdot)).$
\begin{thm}\textbf{(The Necessary Maximum Principle)}
Let $\hat{u}(\cdot)\in {\cal A}$ be a local minimum for
 performance functional
$J(v(\cdot))$(see (2.2)) in the sense that for all bounded
$\beta(\cdot)\in {\cal A}$, there exists $\delta >0$ such that
$\hat{u}(\cdot)+y\beta(\cdot)\in {\cal A}$ for all $y\in
(-\delta,\delta)$ and $h(y):=J(\hat{u}(\cdot)+y\beta(\cdot)),y\in
(-\delta,\delta)$ is minimal at $y=0$.

Suppose there exists a solution
$(\hat{p}_t,\hat{q}_t,\hat{\beta_t}(\cdot),\hat{k}_t)_{t \geq 0}$ of
the adjoint forward-backward stochastic differential   equations
(3.2) corresponding to the admissible pair $(\hat{u}_t;
\hat{x}_t,\hat{y}_t,\hat{z}_t, \hat{r}_t(\cdot))$. Moreover, suppose
that if
 $\hat{X}_t^{1}=X_t^{(\hat{u},~\beta)},\hat{Y}_t^{1}=Y_t^{(\hat{u},~\beta)},\hat{Z}_t^{1}=Z_t^{(\hat{u},~\beta)},$
 (noting (4.1)(4.2)),
 then
 \begin{equation}
 \displaystyle E\displaystyle\int_0^T(X_t^{(1)})
 ^*\hat{q}_t\hat{q}_t^*X_t^{(1)}dt<+\infty,
\end {equation}
\begin{equation}
 \displaystyle E\displaystyle\int_0^T(X_t^{(1)})^*\int_\mathscr{E}
 \hat{\beta}_t(e)\hat{\beta}^*_t(e)\pi(de)X_t^{(1)}dt<+\infty,
\end {equation}
\begin {equation}
  E\displaystyle\int_0^T(Y_t^{(1)})^*( \nabla_zH \nabla_zH^*)
  (t,\hat {x}_t,\hat{y}_t,
  \hat{z}_t, \hat{r}(\cdot),\hat{v}_t,\hat{p}_t,
 \hat{q}_t,\hat{\beta}(\cdot),\hat{k}_t)Y_t^{(1)}dt<+\infty,
\end {equation}
\begin {equation}
  E\displaystyle\int_0^T(Y_t^{(1)})^*\int_\mathscr{E}( \nabla_{r(\cdot)}H \nabla_{r(\cdot)}H^*)
  (t,\hat {x}_t,\hat{y}_t, \hat{z}_t, \hat{r}(\cdot),\hat{v}_t,\hat{p}_t,
 \hat{q}_t,\hat\beta(\cdot),\hat{k}_t)\pi(de)Y_t^{(1)}dt<+\infty,
\end {equation}
\begin {equation}
E\displaystyle\int_0^T\hat{p}_t^*(\xi_t\xi_t^*)dt<+\infty,
\end {equation}
\begin {equation}
E\displaystyle\int_0^T\hat{p}_t^*\int_\mathscr{E}(\zeta_t\zeta_t^*)\pi(d
e)\hat{p}_tdt<+\infty,
\end {equation}
\begin {equation}
  E\displaystyle\displaystyle\int_0^T\hat{k}_t^*(Z_t^{(1)}Z_t^{(1)~*})
  \hat {k}_tdt<+\infty,
\end {equation}
\begin {equation}
  E\displaystyle\displaystyle\int_0^T\hat{k}_t^*
  \int_\mathscr E(R_t^{(v)}(e)R_t^{(v)~*}(e))\pi(de)\hat
  {k}_tdt<+\infty.
\end {equation}

where$$
\begin{array}{ll}
 \xi_t&=\nabla_xg(t,\hat{x}_t,\hat{y}_t,\hat{z}_t,\hat{r}_t(\cdot),\hat{u}_t)X_t^{1}
 +\nabla_yg(t,\hat{x}_t,\hat{y}_t,\hat{z}_t,\hat{r}_t(\cdot),\hat{u}_t)Y_t^{1}
\\\displaystyle
&~~~~~~~+\nabla_zg(t,\hat{x}_t,\hat{y}_t,\hat{z}_t,\hat{r}_t(\cdot),\hat{u}_t)Z_t^{1}
+\displaystyle\int_\mathscr{E}\nabla_{r(\cdot)}g(t,\hat{x}_t,\hat{y}_t,\hat{z}_t,
\hat{r}_t(\cdot),\hat{u}_t)R_t(\bar{e})\pi(d\bar{e}),
\\&~~~~~~~~+\nabla_vg(t,\hat{x}_t,\hat{y}_t,\hat{z}_t,\hat{r}_t(\cdot),\hat{u}_t)\theta_t,\\
\zeta_t&=\nabla_x\sigma(t,\hat{x}_t,\hat{y}_t,\hat{z}_t,\hat{r}_t(\cdot),\hat{u}_t,e)X_t^{1}
+\nabla_y\sigma(t,\hat{x}_t,\hat{y}_t,\hat{z}_t,\hat{r}_t(\cdot),\hat{u}_t,e)Y_t^{1}
\\\displaystyle
&~~~~~~~+\nabla_z\sigma(t,\hat{x}_t,\hat{y}_t,\hat{z}_t,\hat{r}_t(\cdot),\hat{u}_t,
e)Z_t^{1}
+\displaystyle\int_\mathscr{E}\nabla_{r(\cdot)}\sigma(t,\hat{x}_t,\hat{y}_t,\hat{z}_t,
\hat{r}_t(\cdot),\hat{u}_t,e)R_t(\bar{e})\pi(d\bar{e})
\\&~~~~~~~~+\nabla_v\sigma(t,\hat{x}_t,\hat{y}_t,
\hat{z}_t,\hat{r}_t(\cdot),\hat{u}_t,e)\theta_t.
\end{array}$$
 Then $\hat{u}(\cdot)$ is a stationary
point for $E[H|\varepsilon_t]$ in the sense that for almost all
$t\in [0,T]$, we have
$$
E\bigg[H_u(t,\hat{x}_t,\hat{y}_t,\hat{z}_t,\hat{r}_t(\cdot),\hat{u}_t,\hat{p}_t,\hat{q}_t,\hat{\beta}_t(\cdot),
\hat{k}_t)\bigg|\varepsilon_t\bigg]=0.
$$
\end{thm}
\begin{proof}
  From the fact that $h(y)$  is minimal at $y=0$, we
have
\begin{equation}
\begin{array}{ll}
&0=h^{'}(0)=E\displaystyle\int_0^T\bigg(\nabla_xl(t,\hat{x}_t,\hat{y}_t,\hat{z}_t,\hat{r}_t(\cdot),
\hat{u}_t),\hat{X}_t^{1}\bigg)dt\\
&+E\displaystyle\int_0^T\bigg(\nabla_yl(t,\hat{x}_t,\hat{y}_t,\hat{z}_t,\hat{r}_t(\cdot),\hat{u}_t),
\hat{Y}_t^{1}\bigg)dt
+E\displaystyle\int_0^T\bigg(\nabla_zl(t,\hat{x}_t,\hat{y}_t,\hat{z}_t,\hat{r}_t(\cdot),\hat{u}_t),
\hat{Z}_t^{1}\bigg)dt\\&
+E\displaystyle\int_0^T\bigg(\nabla_vl(t,\hat{x}_t,\hat{y}_t,\hat{z}_t,\hat{r}_t(\cdot),\hat{u}_t),
\theta_t\bigg)dt
\\
&+E\displaystyle\int_0^T\displaystyle\int_\mathscr{E}\bigg(\nabla_{r(\cdot)}l(t,\hat{x}_t,\hat{y}_t,
\hat{z}_t,\hat{r}_t(\cdot),\hat{u}_t), \hat
{R}^1_t(\bar{e})\bigg)\pi(d\bar{e})dt
\\&+\displaystyle
E\bigg(\nabla_x\varphi(\hat{x}_T),\hat{X}_T^{1}\bigg)
+E\bigg(h\nabla_y(\hat{y}_0),\hat{Y}_0^{1}\bigg).
\end{array}
\end{equation}
Applying Integration by parts (Lemma 3.1) to
$\big(\hat{p}_t,\hat{X}_t^{1}\big)+\big(\hat{k}_t,\hat{Y}_t^{1}\big)$,
we obtain the following relations
$$
\begin{array}{ll}
&E\bigg(\nabla_x\varphi(\hat{x}_T),\hat{X}_T^{1}\bigg)+E\bigg(\nabla_yh(\hat{y}_0,),\hat{Y}_0^{1}\bigg)
=E\bigg(\hat{p}_T,\hat{X}_T^{1}\bigg)
+E\bigg(-\hat{k}_0,\hat{Y}_0^{1}\bigg)\\
=&-E\displaystyle\int_0^T\bigg(\nabla_xH(t,\hat{x}_t,\hat{y}_t,\hat{z}_t,\hat{r}_t(\cdot),
\hat{u}_t,\hat{p}_t, \hat{q}_t,\hat{\beta}_t(\cdot),\hat{k}_t), \hat{X}_t^{1}\bigg)dt\\
&-E\displaystyle\int_0^T\bigg(\nabla_yH(t,\hat{x}_t,\hat{y}_t,\hat{z}_t,\hat{r}_t(\cdot),\hat{u}_t,\hat{p}_t,
\hat{q}_t,\hat{\beta}_t(\cdot), \hat{k}_t),
\hat{Y}_t^{1}\bigg)dt\\
&-E\displaystyle\int_0^T\bigg(\nabla_zH(t,\hat{x}_t,\hat{y}_t,\hat{z}_t,\hat{r}_t(\cdot),\hat{u}_t,\hat{p}_t,
\hat{q}_t, \hat{\beta}_t(\cdot),\hat{k}_t),\hat{Z}_t^{1}\bigg)dt
\\
&-E\bigg[\displaystyle\int_0^T\displaystyle\int_\mathscr{E}\bigg(\nabla_{r(\cdot)}H(t,\hat{x}_t,
\hat{y}_t,\hat{z}_t,\hat{r}(\cdot),\hat{u}_t,\hat{p}_t,
\hat{q}_t,\hat{\beta}(\cdot),\hat{k}_t),R^1_t(\bar{e})\bigg)\pi(d\bar{e})dt\bigg]\\
&+E\bigg[\displaystyle\int_0^T\bigg(\hat{p}_t,
\nabla_xb(t,\hat{x}_t,\hat{y}_t,\hat{z}_t,\hat{r}_t(\cdot),\hat{u}_t)X_t^{1}+\nabla_yb(t,\hat{x}_t,
\hat{y}_t,\hat{z}_t,\hat{r}_t(\cdot),\hat{u}_t)Y_t^{1}\end {array}
$$
$$
\begin{array}{ll}
&+\nabla_zb(t,\hat{x}_t,\hat{y}_t,\hat{z}_t,\hat{r}_t(\cdot),\hat{u}_t)Z_t^{1}
+\displaystyle\int_\mathscr{E}\nabla_{r(\cdot)}b(t,\hat{x}_t,\hat{y}_t,\hat{z}_t,\hat{r}_t(\cdot),
\hat{u}_t)R_t(\bar{e})\pi(d\bar{e})
\\&+\nabla_vb(t,\hat{x}_t,\hat{y}_t,\hat{z}_t,\hat{r}_t(\cdot),\hat{u}_t)\theta_t\bigg)\bigg]dt
+E\bigg[\displaystyle\int_0^T\bigg(\hat{q}_t,\nabla_xg(t,\hat{x}_t,\hat{y}_t,\hat{z}_t,\hat{r}_t(\cdot),
\hat{u}_t)X_t^{1}\\
&+\nabla_yg(t,\hat{x}_t,\hat{y}_t,\hat{z}_t,\hat{r}_t(\cdot),\hat{u}_t)Y_t^{1}
\displaystyle
+\nabla_zg(t,\hat{x}_t,\hat{y}_t,\hat{z}_t,\hat{r}_t(\cdot),\hat{u}_t)Z_t^{1}
\\
&+\displaystyle\int_\mathscr{E}\nabla_{r(\cdot)}g(t,\hat{x}_t,\hat{y}_t,\hat{z}_t,\hat{r}_t(\cdot),
\hat{u}_t)R_t(\bar{e})\pi(d\bar{e})
+\nabla_vg(t,\hat{x}_t,\hat{y}_t,\hat{z}_t,\hat{r}_t(\cdot),\hat{u}_t)\theta_t\bigg)\bigg]dt
\\
&+E\bigg[\displaystyle\int_0^T\displaystyle\int_\mathscr{E}\bigg(\hat{\beta}_{t}(e),\nabla_x\sigma(t,
\hat{x}_t,\hat{y}_t,\hat{z}_t,\hat{r}_t(\cdot),\hat{u}_t,e)X_t^{1}+\nabla_y\sigma(t,\hat{x}_t,\hat{y}_t,
\hat{z}_t,\hat{r}_t(\cdot),\hat{u}_t,e)Y_t^{1}
\\\displaystyle
&+\nabla_z\sigma(t,\hat{x}_t,\hat{y}_t,\hat{z}_t,\hat{r}_t(\cdot),\hat{u}_t,
e)Z_t^{1}
+\displaystyle\int_\mathscr{E}\nabla_{r(\cdot)}\sigma(t,\hat{x}_t,\hat{y}_t,\hat{z}_t,\hat{r}_t(\cdot),
\hat{u}_t,e)R_t(\bar{e})\pi(d\bar{e})
\\&+\nabla_v\sigma(t,\hat{x}_t,\hat{y}_t,\hat{z}_t,\hat{r}_t(\cdot),\hat{u}_t,e)\theta_t\bigg)\pi(de)dt\bigg]
+E\displaystyle\int_0^T\bigg(\hat{k}_t,-\nabla_xf(t,\hat{x}_t,\hat{y}_t,\hat{z}_t,\hat{r}_t(\cdot),\hat{u}_t)X_t^{1}
\\&-\nabla_yf(t,\hat{x}_t,\hat{y}_t,\hat{z}_t,\hat{r}_t(\cdot),\hat{u}_t)Y_t^{1}
-\nabla_zf(t,\hat{x}_t,\hat{y}_t,\hat{z}_t,\hat{r}_t(\cdot),\hat{u}_t)Z_t^{1}
\\&-\displaystyle\int_\mathscr{E}\nabla_{r(\cdot)}f(t,\hat{x}_t,\hat{y}_t,\hat{z}_t,\hat{r}_t(\cdot),
\hat{u}_t)R_t^{1}(\bar{e})\pi(d\bar{e})
-\nabla_vf(t,\hat{x}_t,\hat{y}_t,\hat{z}_t,\hat{r}_t(\cdot),\hat{u}_t)\theta_t\bigg)dt\\
\end{array}
$$
\begin{equation}
\begin{array}{ll}
=&-E\displaystyle\int_0^T\bigg(\nabla_xl(t,\hat{x}_t,\hat{y}_t,\hat{z}_t,\hat{r}_t(\cdot),\hat{u}_t),
\hat{X}_t^{1}\bigg)dt-E\displaystyle\int_0^T\bigg(\nabla_yl(t,\hat{x}_t,\hat{y}_t,\hat{z}_t,
\hat{r}_t(\cdot),\hat{u}_t),\hat{Y}_t^{1}\bigg)dt
\\&-E\displaystyle\int_0^T\bigg(\nabla_zl(t,\hat{x}_t,\hat{y}_t,\hat{z}_t,\hat{r}_t(\cdot),
\hat{u}_t),\hat{Z}_t^{1}\bigg)dt\\
&
-E\displaystyle\int_0^T\displaystyle\int_\mathscr{E}\bigg(\nabla_{r(\cdot)}l(t,\hat{x}_t,
\hat{y}_t,\hat{z}_t,\hat{r}_t(\cdot),\hat{u}_t), \hat
{R}^1_t(e)\bigg)\pi(de)dt
\\
&+E\displaystyle\int_0^T\bigg(\hat{p}_t,\nabla_vb(t,\hat{x}_t,\hat{y}_t,\hat{z}_t,\hat{r}_t(\cdot),
\hat{u}_t)\theta_t\bigg)dt~~+E\displaystyle\int_0^T\bigg(\hat{q}_t,\nabla_vg(t,\hat{x}_t,\hat{y}_t,
\hat{z}_t,\hat{r}_t(\cdot),\hat{u}_t)\theta_t\bigg)dt
\\
&+E\bigg[\displaystyle\int_0^T\displaystyle\int_\mathscr{E}\bigg(\hat{\beta}_{t}(e),\nabla_v\sigma(t,
\hat{x}_t,\hat{y}_t,\hat{z}_t,\hat{r}_t(\cdot),\hat{u}_t,e)\theta_t\bigg)\pi(de)\\
&+E\displaystyle\int_0^T\bigg(\hat{k}_t,
-\nabla_vf(t,\hat{x}_t,\hat{y}_t,\hat{z}_t,\hat{r}_t(\cdot),
\hat{u}_t)\theta_t\bigg)dt,\\
\end{array}
\end{equation}
where the $L^2$ conditions (4.3-4.10) ensure that the stochastic
integrals with respect to the Brownian motion and  the Poisson
random measure  have zero expectations.

Substituting (4.12) into (4.11), we have
\begin{equation}
E\bigg[\int_0^T\bigg(\nabla_vH(t,\hat{x}_t,\hat{y}_t,\hat{z}_t,\hat{r}_t(\cdot),\hat{u}_t,\hat{p}_t,
\hat{q}_t,\hat{\beta}_t(\cdot),\hat{k}_t),\theta_t\bigg)dt\bigg]=0.
\end{equation}
Fix $t\in [0,T]$ and apply the above to $\theta=(0,
\cdots,\theta_i,0,\cdots,0)$ where
$\theta_i(s)=\alpha_i(\omega)\mathcal {X}_{[t,~~t+r]}{(s)},~s\in
[0,T]$, $t+r\leq T$ and $\alpha_i=\alpha_i(\omega)$ is bounded
$\varepsilon_t$-measurable.

Then it follows from (4.13) that
$$
E\bigg[\int_t^{t+r}\nabla_{
v_i}H(s,\hat{x}_s,\hat{y}_s,\hat{z}_s,\hat{r}_s(\cdot),\hat{u}_s,\hat{p}_s,\hat{q}_s,\hat{\beta}_s(\cdot),\hat{k}_s)
\alpha_ids\bigg]=0.
$$
Differentiating with respect to $r$ at $r=0$ gives
$$
E\bigg[\displaystyle\nabla_{
v_i}H(t,\hat{x}_t,\hat{y}_t,\hat{z}_t,\hat{r}_t(\cdot),\hat{u}_t,\hat{p}_t,\hat{q}_t,\hat{\beta}_t(\cdot),\hat{k}_t)
\alpha_i\bigg]=0.
$$
Since this holds for all bounded $\varepsilon_t$-measurable
$\alpha_i$, we conclude that
$$
E[\nabla_{
v_i}H(t,\hat{x}_t,\hat{y}_t,\hat{z}_t,\hat{r}_t(\cdot),\hat{u}_t,\hat{p}_t,\hat{q}_t,\hat{\beta}_t(\cdot),\hat{k}_t)|\varepsilon_t]=0.
$$
The proof of Theorem 4.1 is completed.
\end{proof}

\end{document}